%% file: HFMMfinal.tex
\documentclass{siamart1116}%
\usepackage{amsmath}
\usepackage{amsfonts}
\usepackage{amssymb}
\usepackage{graphicx}%
\input{ex_shared}
\ifpdf
\hypersetup{
  pdftitle={\TheTitle},
  pdfauthor={\TheAuthors}
}
\fi
\begin{document}

\maketitle

\begin{abstract}
In this paper, we will introduce a new heterogeneous fast multipole method
(H-FMM) for 2-D Helmholtz equation in layered media. To illustrate the main
algorithm ideas, we focus on the case of two and three layers in this work.
The key compression step in the H-FMM is based on a fact that the multipole
expansion for the sources of the free-space Green's function can be used also
to compress the far field of the sources of the layered-media or domain
Green's function, and a similar result exists for the translation operators for the
multipole and local expansions. The mathematical error analysis is shown
rigorously by an image representation of the Sommerfeld spectral form of the
domain Green's function. As a result, in the H-FMM algorithm, both the
``multipole-to-multipole" and ``local-to-local" translation operators are the
same as those in the free-space case, allowing easy adaptation of existing
free-space FMM. All the spatially variant information of the domain Green's
function are collected into the ``multipole-to-local" translations and
therefore the FMM becomes ``heterogeneous". The compressed representation
further reduces the cost of evaluating the domain Green's function when
computing the local direct interactions. Preliminary numerical experiments are
presented to demonstrate the efficiency and accuracy of the algorithm with
much improved performance over some existing methods for inhomogeneous media.
Furthermore, we also show that, due to the equivalence between the complex
line image representation and Sommerfeld integral representation of layered
media Green's function, the new algorithm can be generalized to multi-layered
media with minor modification where details for compression formulas,
translation operators, and bookkeeping strategies will be addressed in a
subsequent paper.

\end{abstract}




\section{Introduction}

To compute the interactions of the electromagnetic or acoustic waves with
objects of complex geometry embedded in the multi-layered media, an attractive
numerical method in the engineering community is to reformulate the frequency
domain Helmholtz equation as a boundary integral equation (BIE) using the
layered-media Green's function where the unknowns are only defined on the
surface of the objects. Unlike the translation-invariant free-space Green's
function for wave scattering in homogeneous media, the layered-media Green's
function incorporates the interface and far-field boundary conditions and
becomes a spatially variant function. In this work, we will refer to the
layered-media Green's function as the domain Green's function or domain kernel
function. Subsequently, the boundary integral equations are discretized using
proper numerical integration techniques, for instance, the trapezoidal rule
with end point corrections~\cite{alpert1999hybrid,kapur1997high} in two
dimensions or the Quadrature by Expansion (QBX) technique in higher
dimensions~\cite{klockner2013quadrature}, resulting in a dense linear system
where the matrix describes how the discretized source and target particles
interact through the domain Green's function. In the numerical solver step, an
important building block is the efficient application of this matrix to a
given vector representing the source contributions, to derive the potential
field due to the domain Green's function interactions at the target particle locations.

There exist several strategies to compute integral operators of the spatially
variant domain Green's function efficiently. For simple geometries, for
instance, the half-space or spheres, one technique is to represent the domain
Green's function contribution as the sum of the free-space Green's function
contributions from both the original source and some image points (the spatial
variant properties are incorporated into the locations of the images). This
approximation allows the direct application of existing free-space fast matrix
vector multiplication algorithms specially designed for the free-space
translation-invariant kernels, including the well-developed fast multipole
method (FMM) packages in
\cite{cho2010wideband,debuhr2016dashmm,gimbutas2015computational,zhang2016rec}%
. Representing works along this direction include the classical Kelvin image
for the half-space problem (or spheres) for the perfect conducting media
\cite{maxwell1881treatise}, so the spatially variant domain Green's function
simply consists of two free-space Coulomb potentials, one from the source
charge and one from its image. In the case of dielectric inhomogeneity, such
as a spherical cavity embedded in a dielectric medium, the reaction field from
the media can also be approximated by a small number of image charges
\cite{cai2007extending}\cite{cai2013}. Unfortunately, for more complex geometry, the image
approximations are extremely hard to derive or non-existent, and for a few
special cases including the multi-layered media, the domain Green's functions
are customarily derived as Sommerfeld integral formulas using integral
transformations. Even when the image approximation formulas are available
(e.g., the two-layered media Helmholtz equation), a large number of images is
usually required. In \cite{o2014efficient}, to approximate the interaction of
$6,400$ particles described by the domain Green's function of the 2D Helmholtz
equation with half-space impedance boundary condition, a total of $1,122,960$
additional images were introduced in a hybrid approach, which combines the
image and Sommerfeld integral representations. Other efforts to speed up the
computation of integral operator for layered media Green's functions include
the inhomogeneuous plane wave method \cite{huchew2000} and the cylindrical wave
decomposition of the Green's function in 3-D and the 2D-FMM \cite{cho2012}.

Another approach is to compress the matrix describing the domain Green's
function interactions directly using the fast direct solvers (FDS)
\cite{greengard2009fast,ho2012fast} or closely related $\mathcal{H}$-matrix
theory \cite{hackbusch1999sparse,hackbusch2000sparse}, where the low-rank
structures of the sub-matrices are derived and processed recursively on a
hierarchical tree structure using purely numerical linear algebra techniques.
However, evaluating the domain Green's functions (entries in the matrix)
involves very expensive computation of the Sommerfeld type integrals, and the
compression stage of the FDS is expensive and memory intensive. It is worth
mentioning that the FMM, FDS, and $\mathcal{H}$-matrix are all hierarchical
algorithms that recursively compress the information in a system to low-rank
or low-dimensional forms and transmit the compressed information non-locally
on a hierarchical tree structure. In this paper, we apply this hierarchical
algorithm design philosophy to multi-layered media domain Green's functions
and present a new hierarchical algorithm for evaluating the spatially variant
domain Green' s function interactions. Our algorithm shares many common
features with FMM and FDS algorithms, especially in the information
transmission patterns on the tree structure: the compressed representations
are transmitted through an upward pass from leaf to parent nodes on a
hierarchical tree structure, collected by interacting nodes and stored as
\textquotedblleft local expansions", and then transmitted to the children
nodes in a downward pass.

We start from the 2-D half-space (two-layered) problem where the domain
Green's function can be explicitly represented with the Sommerfeld integrals
or complex line images. The \textquotedblleft complex line image"
representation intuitively reveals how the compression of the interaction
matrix can be performed analytically on a \textquotedblleft transformed"
matrix which only involves the free-space Green's function, and provides
rigorous error analysis using available analytical results from the classical
free-space FMM. As the compressed representation separates the spatially
variant components and spatially invariant free-space kernels in the domain
Green's function, both the \textquotedblleft multipole-to-multipole" and
\textquotedblleft local-to-local" translations from the existing free-space
Helmholtz FMM algorithm can be easily adapted. Unlike the classical spatially
invariant FMM algorithm, all the spatially variant information are collected
in the \textquotedblleft multipole-to-local" translations and the new
algorithm becomes spatially variant. We refer to this new algorithm as the
\textit{Heterogeneous FMM} due to the heterogeneous nature of the
\textquotedblleft multipole-to-local" translations and the use of the
free-space Green's function and similar translations on the hierarchical tree
structure from the classical FMM. We present the algorithm structure and
demonstrate its accuracy and efficiency by comparing with the hybrid method in
Ref. \cite{o2014efficient} for handling inhomogeneous media. More
interestingly, by relating the results from the complex line image
representation to the equivalent Sommerfeld integral representation, we
demonstrate how the algorithm can be generalized to a particular three-layered
media setting where the complex image representation becomes too complicate to
derive and only the Sommerfeld integral representation via integral
transformations \cite{michalski1990electromagnetic} is available.

This paper is organized as follows. In Sec. \ref{sec:2dhelmholtz}, we present
both the complex line image and Sommerfeld integral representations of the
free-space and domain Green's functions for the 2-D Helmholtz equation with
half-space impedance boundary condition, and present the Sommerfeld integral
representation of the domain Green's function for a particular three-layered
media setting where the analytic image formula is not available. In Sec.
\ref{sec:twolayer_alg}, we focus on the half-space impedance boundary
condition problem, and present the hierarchical algorithm for the efficient
evaluation of the spatially variant domain Green's function interactions. We
will discuss the hierarchical tree structure, compression of the complex image
representations as \textquotedblleft multipole expansions", compression of the
local interactions to allow more efficient evaluations of the integrals,
adaptation of the spatially invariant \textquotedblleft
multipole-to-multipole" and \textquotedblleft local-to-local" translations
from existing free-space FMM, analytical formulas for the heterogeneous
\textquotedblleft multipole-to-local" translations and their efficient
evaluations, and present the algorithm structure and some implementation
details. Preliminary numerical results are presented in Sec.
\ref{sec:numerical_results} to demonstrate the algorithm accuracy and
efficiency. Although it is easier to describe the two-layered media algorithm
and perform the error analysis purely using the \textit{complex line image}
representation in Sec. \ref{sec:twolayer_alg}, we also present the
mathematically equivalent results for the Sommerfeld integral representation.
Comparing these two mathematically equivalent representations, we introduce a
numerical framework solely based on the Sommerfeld integral representation in
Sec. \ref{sec:threelayer_fmm} for a particular three-layered media setting,
which allows direct extension to multi-layered media. Finally, we summarize
our results in Sec. \ref{sec:conclusion}.

\section{2-D Helmholtz Equation in Multi-layered Media}

\label{sec:2dhelmholtz}

We present both the \textit{complex line image} and \textit{Sommerfeld
integral} representations of the domain Green's functions for the 2-D
Helmholtz equation in half-space with impedance boundary condition, and the
\textit{Sommerfeld integral representation} for a particular three-layered
media setting where an analytical expression for the image representation may
not be available.

\subsection{Free-space Green's function}

Consider the 2-D Helmholtz equation in free-space
\[
(\Delta+{k}^{2})u(\mathbf{x})=0
\]
with the Sommerfeld radiation condition at $\infty$
\[
\lim_{r\rightarrow\infty}\sqrt{r}\left(  \frac{\partial}{\partial
r}u(\mathbf{x})-i{k}u(\mathbf{x})\right)  =0,
\]
where $\mathbf{x}=(x,y)$ and $r=||\mathbf{x}||$. Its Green's function is given
by the $0^{th}$ order Hankel function of the first kind as
\begin{equation}
g(\mathbf{x},\mathbf{x}_{0})=\frac{i}{4}H_{0}^{(1)}({k}||\mathbf{x}%
-\mathbf{x}_{0}||) \label{eq:freeGreen}%
\end{equation}
which solves the equation
\begin{equation}
-(\Delta+{k}^{2})g(\mathbf{x},\mathbf{x}_{0})=\delta(\mathbf{x}-\mathbf{x}%
_{0}) \label{eq:freeGreenEqn}%
\end{equation}
with the Sommerfeld radiation condition
\[
\lim_{r\rightarrow\infty}\sqrt{r}\left(  \frac{\partial}{\partial
r}g(\mathbf{x},\mathbf{x}_{0})-i{k}g(\mathbf{x},\mathbf{x}_{0})\right)  =0,
\]
where $\delta$ is the 2-D Dirac delta function, $\mathbf{x}_{0}=(x_{0},y_{0}%
)$, $r=||\mathbf{x}-\mathbf{x}_{0}||$, ${k}$ is the wave number, and
$i=\sqrt{-1}$.

The free-space Green's function can be found in the frequency (spectral)
domain by taking the Fourier transform of Eq.~(\ref{eq:freeGreenEqn}) in the
$x$-direction and solving the resulting ODE in the $y$-direction to give its
spectral representation
\begin{equation}
g(\mathbf{x},\mathbf{x}_{0})=\frac{1}{4\pi}\int_{-\infty}^{\infty}%
\frac{e^{-\sqrt{\lambda^{2}-{k}^{2}}|y-y_{0}|}}{\sqrt{\lambda^{2}-{k}^{2}}%
}e^{i\lambda(x-x_{0})}d\lambda. \label{eq:spectralfree}%
\end{equation}
This representation is often referred to as the Sommerfeld identity, which can
be separated into the propagating and evanescent modes for wave number
variable $|\lambda|<k$ (propagating modes) and $|\lambda|>k$ (evanescent modes
as $|y|\rightarrow\infty$), respectively, to arrive at the following form
after some changes of variables
\begin{align}
g(\mathbf{x},\mathbf{x}_{0})=  &  g(\mathbf{x},\mathbf{x}_{0})_{prop}%
+g(\mathbf{x},\mathbf{x}_{0})_{evan}\nonumber\\
=  &  \frac{i}{4\pi}\int_{0}^{\pi}e^{i{k}\left(|y-y_{0}|\sin{\theta} -  \left(  x-x_{0}\right)
\cos{\theta}\right)  }d\theta\label{eq:evanprop}\\
&  +\frac{1}{4\pi}\int_{0}^{\infty}\frac{e^{-t|y-y_{0}|}}{\sqrt{t^{2}+{k}^{2}%
}}\left(  e^{i\sqrt{t^{2}+{k}^{2}}\left(  x-x_{0}\right)  }+e^{-i\sqrt
{t^{2}+{k}^{2}}\left(  x-x_{0}\right)  }\right)  dt\nonumber
\end{align}
for $|y-y_0|>0$.

The free-space Green's function is commonly used in the potential theory,
where solutions of the Helmholtz equation are represented as combinations of
volume and/or layer potentials defined as the convolution of the Green's
function or its derivatives with certain density functions either over the
volume or surface area of a given object. Theoretical properties of the
free-space Green's function and corresponding potentials are also
well-established in existing literature
\cite{colton2013integral,kress1989linear}, and their efficient evaluations can
be carried out using fast algorithms such as the well-developed wide-band fast
multipole method \cite{cheng2006wideband,cho2010wideband}.

\subsection{Domain Green's function for two-layered media}

In layered media, it is usually possible to derive the spatially variant
domain Green's function analytically either using the method of images
(complex image representation) or applying the integral transforms (e.g.,
Laplace and Fourier transforms) to derive the spectral domain representation
(Sommerfeld integral representation). In this subsection, we focus on the 2-D
half-space Helmholtz equation with the impedance boundary condition
\begin{equation}
\frac{\partial u}{\partial\mathbf{n}}-i\alpha u=0 \label{impedBC}%
\end{equation}
which is imposed on the interface defined by $y=0$, and present the complex
image and Sommerfeld integral representations from existing literature (e.g.,
see \cite{kellogg2012foundations}). \newline

{\noindent\textbf{Complex Image Representation.}} We first present the image
representation of the domain Green's function. The domain Green's function for
the interaction of a point source located at $\mathbf{x}_{0}$ with a target
point at $\mathbf{x}$ is usually decomposed as the sum of the free-space
interaction of the source and target points and contribution from a scattered
field $u_{\mathbf{x}_{0}}^{s}(\mathbf{x})$ which can be explicitly represented
in the two-layered media as complex image contributions of the free-space
kernel as%
\begin{equation}
u_{\mathbf{x}_{0}}^{s}(\mathbf{x})=\int_{0}^{\infty}g(\mathbf{x}%
,\mathbf{x}_{0}^{im}-s\hat{y})\tau(s)ds, \label{scatterField}%
\end{equation}
where $\mathbf{x}_{0}^{im}=(x_{0},-y_{0})$, $\hat{y}=(0,1)$, and $\tau(s)$ is
the complex image charge density distribution. By applying the impedance
boundary condition, the image function $\tau(s)$ can be explicitly found by
(see \cite{o2014efficient})
\begin{equation}
\tau(s)=\delta(s)+\mu(s),\text{ \ }s>0, \label{totalImage}%
\end{equation}
where a point image is indicated by the Dirac delta distribution $\delta(s)$ and a line image $\mu(s)$ is complex and
\begin{equation}
\mu(s)=2i\alpha e^{i\alpha\cdot s}. \label{LineImage}%
\end{equation}
As a result, we have%
\begin{equation}
u_{\mathbf{x}_{0}}^{s}(\mathbf{x})=g(\mathbf{x},\mathbf{x}_{0}^{im})+\int
_{0}^{\infty}g(\mathbf{x},\mathbf{x}_{0}^{im}-s\hat{y})\mu(s)ds,
\label{Sfield}%
\end{equation}
where the first term on the right hand side represents the contribution from
the point-image source, and the second term represents the contributions from
the line-images. Therefore, the domain Green's function $u_{\mathbf{x}_{0}%
}(\mathbf{x})$ for the half-space Helmholtz equation with impedance boundary
condition can be represented in terms of the free-space Green's function
$g(\mathbf{x},\mathbf{x}_{0})$ as
\begin{align}
u_{\mathbf{x}_{0}}(\mathbf{x})  &  =g(\mathbf{x},\mathbf{x}_{0})+u_{\mathbf{x}%
_{0}}^{s}(\mathbf{x})\label{totalfield}\\
&  \equiv g(\mathbf{x},\mathbf{x}_{0})+\left(  g(\mathbf{x},\mathbf{x}%
_{0}^{im})+\int_{0}^{\infty}g(\mathbf{x},\mathbf{x}_{0}^{im}-s\hat{y}%
)\mu(s)ds\right)  .\nonumber
\end{align}

{\noindent\textbf{Sommerfeld Integral Representation.}} The scattered field in
(\ref{Sfield}) involves an integration of an oscillatory line image density
$\mu(s)=2i\alpha e^{i\alpha\cdot s}$, which cannot be handled efficiently with
numerical quadratures directly as in the case for the Laplace equation in
\cite{cai2007extending}. However, using the Sommerfeld identity for
$g(\mathbf{x},\mathbf{x}_{0})$ in Eq.~(\ref{eq:spectralfree}), we can resolve
this difficulty with an analytic integration of the $s$ variable as
follows:
\begin{align}
&  \int_{0}^{\infty}g(\mathbf{x},\mathbf{x}_{0}^{im}-\eta\hat{y}%
)e^{i\alpha\cdot s}ds\nonumber\\
&  =\int_{0}^{\infty}\left[  \frac{1}{4\pi}\int_{-\infty}^{\infty}%
\frac{e^{-\sqrt{\lambda^{2}-{k}^{2}}|y+y_{0}+s|}}{\sqrt{\lambda^{2}-{k}^{2}}%
}e^{i\lambda(x-x_{0})}d\lambda\right]  e^{i\alpha\cdot s}ds\nonumber\\
&  =\frac{1}{4\pi}\int_{-\infty}^{\infty}\frac{e^{-\sqrt{\lambda^{2}-{k}^{2}%
}\left(  y+y_{0}\right)  }e^{i\lambda(x-x_{0})}}{\sqrt{\lambda^{2}-{k}^{2}}%
}\left[  \int_{0}^{\infty}e^{-\sqrt{\lambda^{2}-{k}^{2}}s}e^{i\alpha\cdot
s}ds\right]  d\lambda\nonumber\\
&  =\frac{1}{4\pi}\int_{-\infty}^{\infty}\frac{e^{-\sqrt{\lambda^{2}-{k}^{2}%
}\left(  y+y_{0}\right)  }e^{i\lambda(x-x_{0})}}{\sqrt{\lambda^{2}-{k}^{2}}%
}\frac{1}{\sqrt{\lambda^{2}-{k}^{2}}-i\alpha}d\lambda. \label{ExactInteg}%
\end{align}

Plugging Eqs. (\ref{ExactInteg}) and (\ref{eq:spectralfree}) into Eq.
(\ref{totalfield}), we obtain the following spectral domain representation for
the scattered field (assume $y>0$)%
\begin{equation}
u_{\mathbf{x}_{0}}^{s}(\mathbf{x})=\frac{1}{4\pi}\int_{-\infty}^{\infty}%
\frac{e^{-\sqrt{\lambda^{2}-{k}^{2}}\left(  y+y_{0}\right)  }}{\sqrt
{\lambda^{2}-{k}^{2}}}e^{i\lambda(x-x_{0})}\frac{\sqrt{\lambda^{2}-{k}^{2}%
}+i\alpha}{\sqrt{\lambda^{2}-{k}^{2}}-i\alpha}d\lambda,
\label{eq:spectraldomain}%
\end{equation}
or by defining%
\begin{equation}
\hat{\sigma}(\lambda)= \frac{\sqrt{\lambda^{2}-{k}^{2}}+i\alpha}{\sqrt
{\lambda^{2}-{k}^{2}}-i\alpha},
\end{equation}
we have%
\begin{equation}
\label{eq:multisommer}u_{\mathbf{x}_{0}}^{s}(\mathbf{x})=\frac{1}{4\pi}%
\int_{-\infty}^{\infty}\frac{e^{-\sqrt{\lambda^{2}-{k}^{2}}y}}{\sqrt
{\lambda^{2}-{k}^{2}}}e^{i\lambda x} e^{-\sqrt{\lambda^{2}-{k}^{2}}y_{0}%
}e^{-i\lambda x_{0} } \hat{\sigma}(\lambda)d\lambda
\end{equation}
where $\hat{\sigma}(\lambda)$ is independent of $\mathbf{x}$ and
$\mathbf{x}_{0}$.

\subsection{Domain Green's function for multi-layered media}

For multi-layered media, the explicit forms of the complex image
representations are in general unavailable and the domain Green's functions
are customarily expressed in terms of the Sommerfeld integrals
\cite{chew1995waves, cho2016efficient, chew1999, michalski1990electromagnetic}%
. We leave the detailed derivation of the Sommerfeld integral type domain
Green's functions for different multi-layered media settings to a subsequent
paper. In the following, we cite the results from Ref. \cite{lai2014fast} for
a particular three-layered media setting, and compare its Sommerfeld integral
representation of the domain Green's function with that of the two-layered
media in Eq.~(\ref{eq:multisommer}).

Assume a point source is located at $\mathbf{x}_{0}=(x_{0},y_{0})$ in the top
layer with a wave number $k_{1}$ and the two interfaces of the three-layered
media are located at $y=0$ and $y=-d$, respectively. The Sommerfeld integral
representation for the scattered field in the top layer ($y>0$) can be
represented as
\begin{equation}
u_{1}^{s}(\mathbf{x})=\frac{1}{4\pi}\int_{-\infty}^{\infty}\frac
{e^{-\sqrt{\lambda^{2}-{k^{2}_{1}}}y}}{\sqrt{\lambda^{2}-{k^{2}_{1}}}%
}e^{i\lambda x}e^{-\sqrt{\lambda^{2}-{k^{2}_{1}}}y_{0}}e^{-i\lambda x_{0}%
}\sigma_{1}(\lambda)d\lambda\label{eq:multitop}%
\end{equation}
where the unknown function $\sigma_{1}(\lambda)$ will be determined
later. In the middle layer with a wave number $k_{2}$, the scattered field
$u_{2}^{s}$ can be written as the sum of the contributions $u_{2}^{t}$ (from
upper interface) and $u_{2}^{b}$ (from lower interface) as
\begin{align}
u_{2}^{t}(\mathbf{x})  &  =\frac{1}{4\pi}\int_{-\infty}^{\infty}\frac
{e^{\sqrt{\lambda^{2}-{k^{2}_{2}}}y}}{\sqrt{\lambda^{2}-{k^{2}_{2}}}%
}e^{i\lambda x}e^{-\sqrt{\lambda^{2}-{k^{2}_{2}}}y_{0}}e^{-i\lambda x_{0}%
}\sigma_{2}^{+}(\lambda)d\lambda,\\
u_{2}^{b}(\mathbf{x})  &  =\frac{1}{4\pi}\int_{-\infty}^{\infty}%
\frac{e^{-\sqrt{\lambda^{2}-{k^{2}_{2}}}(y+d)}}{\sqrt{\lambda^{2}-{k^{2}_{2}}%
}}e^{i\lambda x}e^{-\sqrt{\lambda^{2}-{k^{2}_{2}}}y_{0}}e^{-i\lambda x_{0}%
}\sigma_{2}^{-}(\lambda)d\lambda,
\end{align}
and in the bottom layer with wave number $k_{3}$, we have
\begin{equation}
u_{3}^{s}(\mathbf{x})=\frac{1}{4\pi}\int_{-\infty}^{\infty}\frac
{e^{\sqrt{\lambda^{2}-{k^{2}_{3}}}(y+d)}}{\sqrt{\lambda^{2}-{k^{2}_{3}}}%
}e^{i\lambda x}e^{-\sqrt{\lambda^{2}-{k^{2}_{3}}}y_{0}}e^{-i\lambda x_{0}%
}\sigma_{3}(\lambda)d\lambda,
\end{equation}
where $\sigma_{2}^{+}(\lambda)$, $\sigma_{2}^{-}(\lambda)$, and $\sigma
_{3}(\lambda)$ are unknown quantities which are associated with layer reflection
and transmission coefficients of waves in spectral domain. When the interface conditions are given
by $[u]=0$ and $[\frac{\partial u}{\partial n}]=0$, these quantities can
be determined by solving the linear system
\begin{equation}
\left(
\begin{array}
[c]{cccc}%
\frac{1}{\sqrt{\lambda^{2}-{k^{2}_{1}}}} & \frac{-1}{\sqrt{\lambda^{2}%
-{k^{2}_{2}}}} & -\frac{e^{-\sqrt{\lambda^{2}-{k^{2}_{2}}}d}}{\sqrt
{\lambda^{2}-{k_{2}}^{2}}} & 0\\
0 & \frac{e^{-\sqrt{\lambda^{2}-{k^{2}_{2}}}d}}{\sqrt{\lambda^{2}-{k^{2}_{2}}%
}} & \frac{1}{\sqrt{\lambda^{2}-{k^{2}_{2}}}} & \frac{-1}{\sqrt{\lambda
^{2}-{k^{2}_{3}}}}\\
1 & 1 & -e^{-\sqrt{\lambda^{2}-{k^{2}_{2}}}d} & 0\\
0 & e^{-\sqrt{\lambda^{2}-{k^{2}_{2}}}d} & -1 & -1
\end{array}
\right)  \left(
\begin{array}
[c]{c}%
\sigma_{1}(\lambda)\\
\sigma_{2}^{+}(\lambda)\\
\sigma_{2}^{-}(\lambda)\\
\sigma_{3}(\lambda)
\end{array}
\right)  =\left(
\begin{array}
[c]{c}%
\frac{-1}{\sqrt{\lambda^{2}-{k^{2}_{1}}}}\\
0\\
1\\
0
\end{array}
\right)  .\nonumber
\end{equation}
When all the source and target points are located in the top layer, we see
that the domain Green's function $u_{1}^{s}$ is similar to
Eq.~(\ref{eq:multisommer}) but with a different $\sigma$ function, the nature
of which will be revealed later.

One additional complexity of the multi-layered media computation is that the
source and target points may be located in different layers. We leave the
technical details for different domain Green's functions and their
compressions, translations, and bookkeeping strategies for such cases in a
subsequent paper, and only focus on the Green's function in the form of
Eq.~(\ref{eq:multitop}) in this paper. This simplifies the notations while
presenting the main idea of the new heterogeneous FMM algorithm, and still
covers important applications for scattering of multiple objects over
conducting surface (the two-layered media with impedance boundary condition)
and layered media (meta-surfaces over layered substrate in material sciences).

\subsection{Domain Green's function in integral equation methods}

In most integral equation formulations of the Helmholtz equation, unlike the
translation invariant free-space Green's function that only depends on the
distance of $\mathbf{x}$ and $\mathbf{x}_{0}$, the domain Green's function for
general complex geometry brings complications for being a two variable
function and its values are simply no longer translation invariant but
spatially variant. As a result, the computation and evaluation of the domain
Green's function are more expensive than finding the solution of the original
differential equation. There are a few exceptions, including the simulation of
the layered-media Helmholtz equation where the interface or boundary is
infinite but regular. For such cases, if using the free-space Green's
function, the resulting potentials will involve the evaluations of integrals
(potentials) over infinite interfaces. However, as the geometry is regular,
one can analytically derive the spatially variant domain Green's function, in
the form of a Sommerfeld integral representation using integral
transformations such as the Laplace and Fourier transforms.

There are many advantages by using the domain Green's function in forming the
integral equation method (IEM) for the multi-layered media problem, for
instance, the interface conditions are naturally enforced by the domain
Green's function and no unknowns are necessary on the layer interfaces.
However, the numerical solution of the integral equation poses many challenges
and is still an active research topic. In addition to problems common to all
integral equation approaches such as the design of high order quadrature and
derivation of well-conditioned systems, the IEM for layered media using the
domain Green's function has its specific challenges. In particular, the
evaluation of the domain Green's function interactions with large number of
source and target points is expensive for either the complex image
representation, or the Sommerfeld integral, or even the optimized hybrid
representations. This implies that explicitly constructing the discretized
interaction matrix is also expensive, and therefore matrix compression using
the FDS will be costly where purely numerical linear algebra techniques are
applied. This paper focuses on the fast application of the domain Green's
function to a given density function $\rho(\mathbf{x}_{0})$ as in%
\begin{equation}
\phi(\mathbf{x})=\int u_{\mathbf{x}_{0}}(\mathbf{x)\rho(x}_{0})d\mathbf{x}%
_{0},
\end{equation}
where the integral either represents a volume potential or a surface layer
potential. After discretization, the resulting linear algebra question becomes
how to efficiently calculate the matrix-vector multiplication of $A\mathbf{v}$
where entries in the matrix $A$ are given by the domain Green's function
$A_{i,j}=u_{\mathbf{x}_{j}}(\mathbf{x}_{i})$. The main results of this paper
include (a) the analysis-based low-rank compression of the matrix $A$, which
is not directly performed on the matrix itself as in the FDS methods, but on a
closely related matrix after certain transformations; (b) how the compressed
representations can be transmitted through the hierarchical tree structure
using analysis-based translation operators, and whenever possible, utilizing
existing translation operators for the free-space kernels; and (c) the
selected compression and translation strategies allow the implementation of a
\textquotedblleft heterogeneous" FMM algorithm for the layered media by an
easy adaptation of existing free-space FMM codes.

\section{Algorithm for Two-layered Media}

\label{sec:twolayer_alg}

In this section, we present the technical details of a fast hierarchical
algorithm for the two-layered media domain Green's function. The algorithm is
similar in structure to that of FMM and the compression stage of FDS, and is
developed by considering the design philosophy of the hierarchical modeling
technique. This technique identifies any low-rank, or low-dimensional, or
other compact features in a given system, recursively collects the compressed
representations from children to parents, and transmits the information
between different nodes on a hierarchical tree structure using properly
compressed translation operators. It is worth mentioning that the resulting
hierarchical models are often re-expressed as recursive algorithms, which can
be easily interfaced with existing dynamical schedulers from High-Performance
Computing (HPC) community for optimal parallel efficiency
\cite{HPX5web,blumofe1996cilk,frigo1998implementation}.

In addition to FMM and FDS, different aspects of the hierarchical modeling
technique have been known and addressed by many researchers previously.
Examples include the classical fast Fourier transform (FFT)
\cite{cooley1965algorithm} where the Halving Lemma shows how data can be
compressed and the odd-even term splitting of the polynomials creates a
hierarchical tree to allow recursively processing the compressed information
efficiently; the multigrid method (MG)
\cite{brandt1977multi,hackbusch2013multi} where the hierarchical tree
structure is formed via adaptively refining the computational domain, and data
compression and transmission are performed using the relaxation (smoother) and
projection (restriction) operators by analyzing the frequency domain behaviors
of the error functions between different levels of the (adaptive) tree to
effectively reduce the high frequency errors. When there are $n$ terms (FFT)
in the polynomial or $n$ approximately uniformly distributed particles (MG or
FMM), the depth of the hierarchical tree is normally $\mathcal{O}(\log n)$ and
the number of tree nodes is approximately $\mathcal{O}(n)$. Therefore, if each
level only requires $\mathcal{O}(n)$ operations (e.g., FFT), the algorithm
complexity will be $\mathcal{O}(n\log n)$. If each tree node only requires a
constant amount of operations (e.g. MG or FMM), the algorithm complexity will
be asymptotically optimal $\mathcal{O}(n)$. In this section, we describe our
algorithm following the design guidelines of the hierarchical modeling technique.

\subsection{Adaptive hierarchical tree structure}

Consider the Helmholtz equation in 2-D with half-space impedance boundary
condition for the scattering of a finite-sized object with a complex geometry
in the upper half plane $y>0$. A surface integral equation can be derived to
give the scattering solution as layer potentials through a convolution of the
domain Green's function or its derivatives with some unknown density functions
over the object's surface. We assume the surface is \textquotedblleft
discretized" into a number of particles via proper numerical integration
techniques. In the hierarchical modeling technique, a spatial adaptive
hierarchical tree is first generated. In our algorithm, the tree structure is
identical to that in FMM or FDS and is generated by a recursive partition to
divide the particle-occupant region into nested cubical boxes, where the root
box is the smallest bounding box that contains the entire particle set.
Without loss of generality, the root box is normalized to size $1$ along each
side. The root box is partitioned equally along each dimension. The partition
continues recursively on the resulting box until the box contains no more than
$s$ beads, at which point it becomes a leaf node. Empty boxes encountered
during partition are pruned off. In our implementation, the value $s$ is
chosen depending on the size of the particle set and other factors to allow an
optimal performance.

\subsection{Low rank compression}

\label{sec:lowrank} 
\begin{figure}[t]
\centering  \includegraphics[width=3.2in]{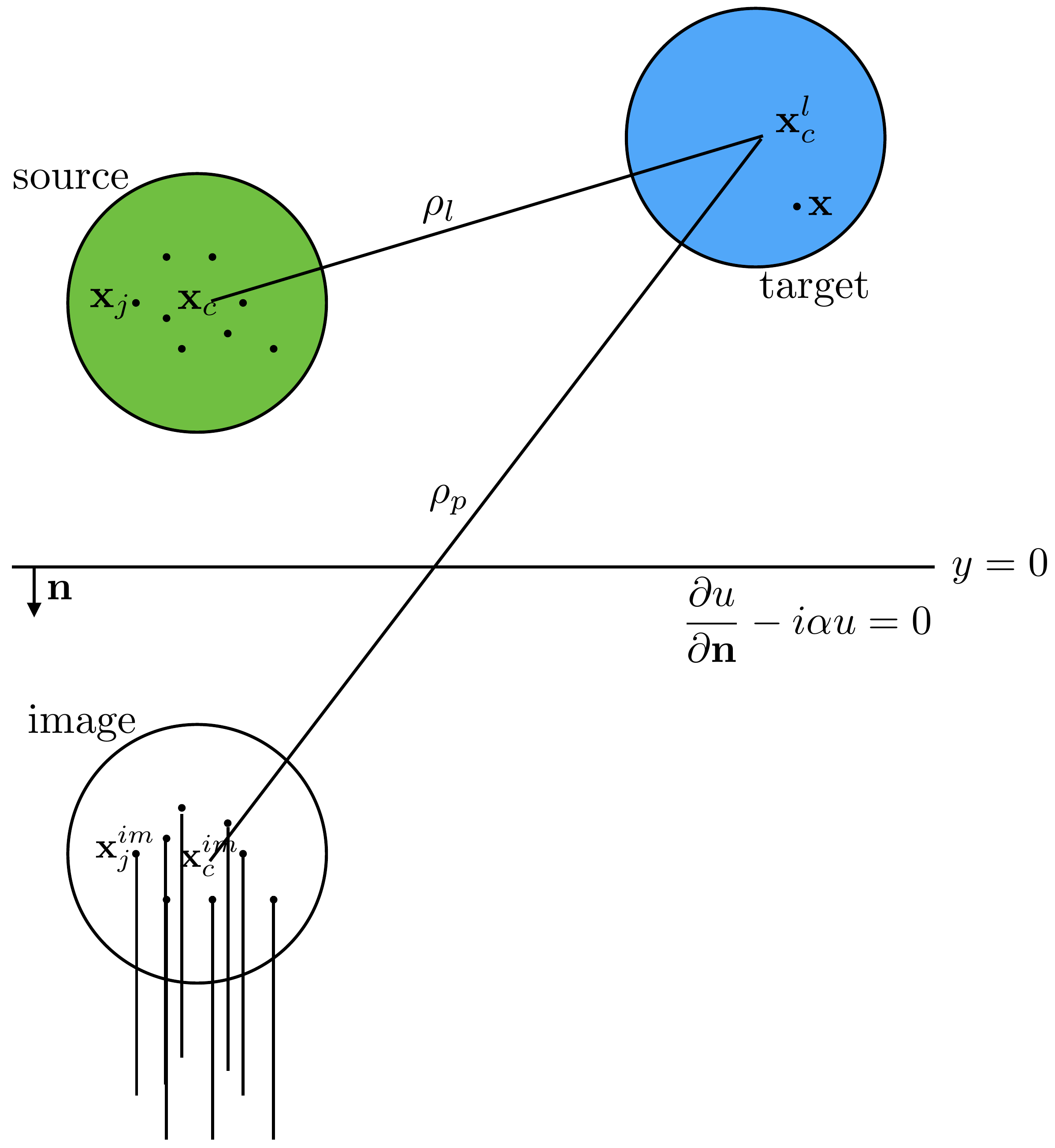} \caption{Impedance
half-space and notation.}%
\label{fig:impedance_half}%
\end{figure}

The low-rank structure for well-separated source and
target points has been extensively studied for the free-space FMM and FDS
algorithms. Consider $N$ sources with strength $q_{j}$ placed at
$\mathbf{x}_{j}=(x_{j},y_{j})$ in a circle centered at $\mathbf{x}_{c}%
=(x_{c},y_{c})$ with radius $R$ on top of the half-space, and suppose we are
interested in the field at $\mathbf{x}$ due to all the source points given by
$u(\mathbf{x})=\sum_{j=1}^{N}q_{j}g(\mathbf{x},\mathbf{x}_{j})$ where
$g(\mathbf{x},\mathbf{x}_{j})$ is the free-space Green's function
contribution. We say $\mathbf{x}$ is well-separated from the sources if the
distance between $\mathbf{x}$ and the source circle center $\mathbf{x}_{c}$ is
at least $3R$, see Fig.~\ref{fig:impedance_half}. \newline

{\noindent\textbf{Free-space Green's function compression.}} Using Graf's
addition theorem \cite{abramowitz1966handbook}, the free-space Green's
function interaction of well-separated sources $\mathbf{x}_{j}$ and target
$\mathbf{x}$ can be compressed as a \textquotedblleft multipole expansion"
given by
\begin{equation}
u(\mathbf{x})\approx\frac{i}{4}\sum_{p=-P}^{P}\alpha_{p}H_{p}({k}%
|\mathbf{x}-\mathbf{x}_{c}|)e^{ip\theta_{c}}, \label{multipole_free}%
\end{equation}
where
\begin{equation}
\alpha_{p}=\sum_{j=1}^{N}q_{j}e^{-ip\theta_{j}}J_{p}({{k}}\rho_{j}),
\label{multi_coeff}%
\end{equation}
$\theta_{c}$ is the polar angle of $\mathbf{x}-\mathbf{x}_{c}$, $(\rho
_{j},\theta_{j})$ are the polar coordinates of the complex number
$\mathbf{x}_{j}-\mathbf{x}_{c}$, and the number of terms $P$ is a constant
independent of the number of the sources $N$. \newline

{\noindent\textbf{Domain Green's function compression: complex image
representation.}} For the half-space problem with an impedance boundary
condition, the field at $\mathbf{x}$ due to all the source points is
\begin{equation}
u(\mathbf{x})=\sum_{j=1}^{N}q_{j}u_{\mathbf{x}_{j}}(\mathbf{x})=\sum_{j=1}%
^{N}q_{j}g(\mathbf{x},\mathbf{x}_{j})+\sum_{j=1}^{N}q_{j}u_{\mathbf{x}_{j}%
}^{s}(\mathbf{x}).
\end{equation}
As the first term on the right hand side represents the free-space Green's
function interaction and is already compressed in Eq.~(\ref{multipole_free}),
we focus on the compression of the second term representing the scattered
field. When using the complex image representation, it can be compressed
simply as follows:%
\begin{align}
& u^{s}(\mathbf{x})=\sum_{j=1}^{N}q_{j}u_{\mathbf{x}_{j}}^{s}(\mathbf{x}%
)\nonumber\\
& =\frac{i}{4}\sum_{j=1}^{N}q_{j}\left(  H_{0}^{(1)}({k}|\mathbf{x}%
-\mathbf{x}_{j}^{im}|)+\int_{0}^{\infty}H_{0}^{(1)}({k}|\mathbf{x}%
-(\mathbf{x}_{j}^{im}-s\hat{y})|)\mu(s)ds\right) \nonumber\\
& \approx\frac{i}{4}\sum_{p=-P}^{P}\bar{\alpha}_{p}\left(  H_{p}%
({k}|\mathbf{x}-\mathbf{x}_{c}^{im}|)e^{ip\theta_{im}}+\int_{0}^{\infty}%
H_{p}({k}|\mathbf{x}-(\mathbf{x}_{c}^{im}-s\hat{y})|)e^{ip\hat{\theta}_{im}%
}\mu(s)ds\right)  \label{eq:greenimage}%
\end{align}
where $\mathbf{x}_{j}^{im}=(x_{j},-y_{j})$ are the coordinates of the
point-image charge, $\bar{\alpha}_{p}$ is the complex conjugate of the
free-space multipole coefficient $\alpha_{p}$ in Eq.~(\ref{multi_coeff}) (see
\cite{o2014efficient}), and%
\begin{align}
&\theta_{im} \text{ is the polar angle of complex number  }\mathbf{x}%
-\mathbf{x}_{c}^{im},\\
&\hat{\theta}_{im} \text{ is the polar angle of complex number }\mathbf{x}%
-(\mathbf{x}_{c}^{im}-s\hat{y}).
\end{align}

Therefore, the \textquotedblleft multipole expansion" for both the original
and image sources is
\begin{align}
u(\mathbf{x})\approx &  \frac{i}{4}\sum_{p=-P}^{P}\alpha_{p}H_{p}%
({k}|\mathbf{x}-\mathbf{x}_{c}|)e^{ip\theta_{c}}+\frac{i}{4}\sum_{p=-P}%
^{P}\bar{\alpha}_{p}\left(  H_{p}({k}|\mathbf{x}-\mathbf{x}_{c}^{im}%
|)e^{ip\theta_{im}}\right. \nonumber\\
&  +\int_{0}^{\infty}\left.  H_{p}({k}|\mathbf{x}-(\mathbf{x}_{c}^{im}%
-s\hat{y})|)e^{ip\hat{\theta}_{im}}\mu(s)ds\right)  , \label{greens_function}%
\end{align}
which is the key formula behind the proposed new heterogeneous FMM for layered media.

We emphasize that the number of terms $P$ for the scattered field expansion is
the same as the one in the free-space expansion for the same accuracy
requirement. This can be rigorously justified by the observation that when the
original sources $\mathbf{x}_{j}$ (in the green circle centered at
$\mathbf{x}_{c}$) in Fig.~\ref{fig:impedance_half} are well-separated from the
target point $\mathbf{x}$ (in the blue circle centered at $\mathbf{x}_{c}^{l}%
$), all the corresponding point-images $\mathbf{x}_{j}^{im}$ (in the circle
centered at $\mathbf{x}_{c}^{im}$) are also well-separated from $\mathbf{x}$,
and this is true also for the set of line-images on the rays starting from
$\mathbf{x}_{c}^{im}$ with the same $s$ value. \newline

{\noindent\textbf{Remark 1:}} Eq.~(\ref{greens_function}) suggests that for
the domain Green's function interactions, when the source and target clusters
are well-separated, it is possible to only compress the translation invariant
free-space Green's function using a $P$-term multipole expansion with
coefficients $\alpha_{p}$ as in Eq.~(\ref{multipole_free}) for a prescribed
accuracy requirement, and all other related information in
Eq.~(\ref{greens_function}) can be recovered from $\alpha_{p}$ to the same
accuracy. Also, unlike the FDS, the compression of the domain Green's function
is not performed directly on the matrix entries, but on another matrix after
some spatially variant transformations implicitly described in
Eq.~(\ref{greens_function}), and these transformations involve the complex
conjugation operator. Finally, as the compression is only on the free-space
Green's function, deriving parent's compressed representation and
corresponding error analysis are exactly the same as those in the classical
FMM algorithms, where the \textquotedblleft multipole-to-multipole"
translation can be used without any modification. We therefore skip the
details of this translation operator in this paper. \newline

{\noindent\textbf{Domain Green's function compression: Sommerfeld integral
representation.}} For multi-layered media, deriving the image representation
of the domain Green's function becomes very complicated, and most existing
techniques apply the integral transformations to obtain the Sommerfeld
representation directly. Instead of deriving the compression formula using the
Sommerfeld integral domain Green's function in Eq.~(\ref{eq:spectraldomain}),
we reformulate its compressed image representation in
Eq.~(\ref{greens_function}) as an equivalent Sommerfeld integral
representation for the sources of $N$ particles. Then, it is compared directly
with the uncompressed Sommerfeld representations. This comparison provides
insight on a direct compression of the domain Green's function in the spectral domain.

We start from the Sommerfeld representation of $H_{n}({{k}}\rho)e^{in\theta}$.
Applying the relation
\[
H_{n}({k}\rho)e^{in\theta}=\left(  -\frac{1}{{{k}}}\right)  ^{n} \left(
\frac{\partial}{\partial x}+i\frac{\partial}{\partial y}\right)  ^{n}H_{0}%
({k}\rho)
\]
where $(\rho,\theta)$ are the polar coordinates of the complex number $x+iy$,
and using the Sommerfeld representation of $H_{0}({k}\rho)$ given in
Eq.~(\ref{eq:spectralfree}), we have for $y>0$,
\begin{equation}
H_{n}({k}\rho)e^{in\theta}=  \frac{(-i)^n} {i \pi}         \int_{-\infty}^{\infty
}\frac{e^{-\sqrt{\lambda^{2}-{k}^{2}}y}}{\sqrt{\lambda^{2}-{k}^{2}}%
}e^{i\lambda x}\left(  \frac{\lambda-\sqrt{\lambda^{2}-{k}^{2}}}{{k}}\right)
^{n}d\lambda. \label{eq:higherordersommer}%
\end{equation}
Plugging this representation into Eq.~(\ref{greens_function}), and integrating
the $s$ variable analytically, we have the compressed Sommerfeld
representation directly as
\begin{align}
u^{s}(\mathbf{x})\approx &  \int_{-\infty}^{\infty}\underbrace{\frac
{e^{-\sqrt{\lambda^{2}-{k}^{2}}(y+y_{c})}}{\sqrt{\lambda^{2}-{k}^{2}}%
}e^{i\lambda(x-x_{c})}}_{\mbox{free-space info}}\underbrace{\left(  \frac
{1}{4\pi}\sum_{p=-P}^{P}\bar{\alpha}_{p}(-i)^{p}\left(  \frac{\lambda
-\sqrt{\lambda^{2}-{k}^{2}}}{{k}}\right)  ^{p}\right)  }_{\mbox{compressed}}%
\label{eq:compressed}\\
&  \underbrace{\left(  \frac{\sqrt{\lambda^{2}-{k}^{2}}+i\alpha}{\sqrt
{\lambda^{2}-{k}^{2}}-i\alpha}\right)  }_{\mbox{image info}}d\lambda.\nonumber
\end{align}
Comparing with the uncompressed Sommerfeld representation by adding up the
Sommerfeld representation of the scattered field in
Eq.~(\ref{eq:spectraldomain}) for each source in $u^{s}(\mathbf{x})=\sum
_{j=1}^{N}q_{j}u_{\mathbf{x}_{j}}^{s}(\mathbf{x})$, we have
\begin{align}
u^{s}(\mathbf{x})  &  =\int_{-\infty}^{\infty}\underbrace{\frac{e^{-\sqrt
{\lambda^{2}-{k}^{2}}(y+y_{c})}}{\sqrt{\lambda^{2}-{k}^{2}}}e^{i\lambda
(x-x_{c})}}_{\mbox{free-space info}}\underbrace{\left(  \frac{1}{4\pi}\sum
_{j=1}^{N}q_{j}e^{-\sqrt{\lambda^{2}-{k}^{2}}(y_{j}-y_{c})}e^{i\lambda
(x_{c}-x_{j})}\right)  }_{\mbox{uncompressed}}\label{eq:uncompressed}\\
&  \underbrace{\left(  \frac{\sqrt{\lambda^{2}-{k}^{2}}+i\alpha}{\sqrt
{\lambda^{2}-{k}^{2}}-i\alpha}\right)  }_{\mbox{image info}}d\lambda.\nonumber
\end{align}
We further notice that $\bar{\alpha}_{p}$ is independent of $\lambda$, and the
compressed term in Eq.~(\ref{eq:compressed}) is the Laurent expansion in
$z=\frac{\lambda-\sqrt{\lambda^{2}-{k}^{2}}}{{k}}$ of the uncompressed term in
Eq.~(\ref{eq:uncompressed}).\newline


{\noindent\textbf{Remark 2:}} Comparing Eq.~(\ref{eq:compressed}) with
Eq.~(\ref{eq:uncompressed}), we can identify the roles of different terms in
the Sommerfeld representation of the domain Green's function. In particular,
we see that the conjugates of the free-space multipole expansion coefficients
are the same as the Laurent expansion ones in the compressed representation.
This observation reveals how the domain Green's function can be compressed
directly when the line image $\mu(s)$ in Eq.~(\ref{LineImage}) is
complicated or unavailable. For instance, for the top layer Sommerfeld domain
Green's function in Eq.~(\ref{eq:multitop}) for the three layered media
setting, as the terms \textquotedblleft free-space info" and \textquotedblleft
uncompressed" have the same structure as in Eq.~(\ref{eq:uncompressed}), and
the term \textquotedblleft image info" is independent of $\mathbf{x}$ and
$\mathbf{x}_{0}$, we can therefore simply compute the free-space multipole
expansion either directly from the sources, or through the free-space
\textquotedblleft multipole-to-multipole" translations, and the results will
directly give a compressed Sommerfeld representation similar to
Eq.~(\ref{eq:compressed}). In Sec. \ref{sec:threelayer_fmm}, we apply this
observation and present an algorithm framework purely based on the Sommerfeld representation.

\medskip
{\noindent\textbf{Remark 3:}} The error analysis of the direct compression of
the Sommerfeld integral representation described in Eqs.~(\ref{eq:compressed})
and (\ref{eq:uncompressed}) is not an easy task. Luckily for the 2-D
half-space problem, the error analysis becomes trivial when performed in the
physical domain using the image representation. Note that the error analysis
only requires that the target box and all the image sources are
well-separated, which can be easily carried out for the multi-layered case in
Eq.~(\ref{eq:multitop}) using repeated image reflections, without knowing the
exact density and location of the image. \newline

{\noindent\textbf{Local expansions for received information.}} We have so far
discussed how the domain Green's function can be compressed in different ways
so the compressed representations can be transmitted to the ``receiving" nodes
on the hierarchical tree structure. In the hierarchical algorithms, the
received information is also stored in some compressed compact form. For the
Helmholtz equation, the most convenient and commonly used form is the Bessel
functions based ``local expansion" as in the classical FMM algorithms. We
adopt this representation in this paper to store the received information, so
the ``local-to-local" translations in the free-space FMM algorithms can be
used without any modification. Other compact forms are also possible, for
example, the ``equivalent source" representation in the kernel independent FMM
algorithms \cite{ying2004kernel,ying2003new}.

\subsection{Translations on the hierarchical tree structure}

We discuss how the compressed representations can be transmitted on the
hierarchical tree structure in this section. As our selected ``multipole" and
``local" representations of the compressed domain Green's function are the
same as those for the translation invariant free-space Green's function,
existing ``multipole-to-multipole" and ``local-to-local" translations in the
free-space FMM can be applied without any modification. We therefore focus on
the ``multipole-to-local" (M2L) translation operator, and study how the
multipole expansion of the compressed domain Green's function can be converted
to local expansions.

We start from the following well-known M2L translation operator for the
free-space Green's function. Consider the same source points $\mathbf{x}_{j}$,
$j=1,\cdots,N$ described in Fig.~\ref{fig:impedance_half} and the compressed
representation of the free-space kernel in Eq.~(\ref{multipole_free}). Then
the potential $u(\mathbf{x})$ can be translated to a local expansion using
Graf's addition theorem as
\begin{equation}
u(\mathbf{x})\approx\frac{i}{4}\sum_{p=-P}^{P}\beta_{p}^{f}J_{p}%
({k}|\mathbf{x}-\mathbf{x}_{c}^{l}|)e^{i{k}\theta},
\end{equation}
where the coefficients are
\begin{equation}
\beta_{p}^{f}=\sum_{m=-P}^{P}\alpha_{m}H_{m-p}({k}\rho_{l})e^{i(m-p)\theta
_{l}},
\end{equation}
$\theta$ is the polar angle of $\mathbf{x}-\mathbf{x}_{c}^{l}$, and $(\rho
_{l},\theta_{l})$ are the polar coordinates of $\mathbf{x}_{c}^{l}%
-\mathbf{x}_{c}$. Because the complex image representation of the domain
Green's function is given in terms of the free-space Green's function, we can
therefore plug the free-space M2L translation formula in the compressed image
representation of the scattered field
\begin{align}
u^{s}(\mathbf{x})  &  \approx\frac{i}{4}\sum_{p=-P}^{P}\bar{\alpha}_{p}\bigg(
H_{p}({k}|\mathbf{x}-\mathbf{x}_{c}^{im}|)e^{ip\theta_{im}}\bigg. \\
&  \left.  +\int_{0}^{\infty}H_{p}({k}|\mathbf{x}-(\mathbf{x}_{c}^{im}%
-s\hat{y})|)\mu(s)e^{ip\hat{\theta}_{im}}ds\right)  ,\nonumber
\end{align}
to derive its local expansion given by
\begin{align}
u^{s}(\mathbf{x})  &  =\frac{i}{4}\sum_{p=-P}^{P}\sum_{m=-p}^{p}\bar{\alpha
}_{m}\bigg(  H_{m-p}({k}\tilde{\rho}_{im})J_{p}({k}|\mathbf{x}-\mathbf{x}%
_{c}^{l}|)e^{i(m-p)\tilde{\theta}_{im}}e^{ip\theta}\bigg.. \nonumber\\
&  \left.  +\int_{0}^{\infty}H_{m-p}({k}\hat{\tilde{\rho}}_{im})J_{p}%
({k}|\mathbf{x}-\mathbf{x}_{c}^{l}|)e^{i(m-p)\hat{\tilde{\theta}}_{im}}%
\mu(s)e^{ip\theta}ds\right) \nonumber\\
&  =\frac{i}{4}\sum_{p=-P}^{P}\beta_{p}^{s}J_{p}({k}|\mathbf{x}-\mathbf{x}%
_{c}^{l}|)e^{ip\theta},
\end{align}
where the local expansion coefficients are given by
\begin{equation}
\beta_{p}^{s}=\sum_{m=-p}^{p}\bar{\alpha}_{m}\left(  H_{m-p}({k}\tilde{\rho
}_{im})e^{i(m-p)\tilde{\theta}_{im}}+\int_{0}^{\infty}H_{m-p}({k}\hat
{\tilde{\rho}}_{im})e^{i(m-p)\hat{\tilde{\theta}}_{im}}\mu(s)ds\right)  ,
\end{equation}
$\theta$ is the polar angle of $\mathbf{x}-\mathbf{x}_{c}^{l}$, and
\begin{align}
&(\tilde{\rho}_{im},\tilde{\theta}_{im}) \text{ are the polar coordinates of }
\mathbf{x}_{c}^{l}-\mathbf{x}_{c}^{im},\\
&(\hat{\tilde{\rho}}_{im},\hat{\tilde{\theta}}_{im}) \text{ are the polar
coordinates of } \mathbf{x}_{c}^{l}-(\mathbf{x}_{c}^{im}-s\hat{y}).
\end{align}
%
The local expansion for $u(\mathbf{x})$ is simply the sum of the free-space
Green's function and scattered field local expansions. As the translation
operator from the compressed \textquotedblleft multipole coefficients"
$\{\alpha_{p}\}$ to the local coefficients $\{\beta_{p}^{s}\}$ involves the
complex conjugate operator, for notation reasons, instead of combining the
free-space with the complex image contributions in one single translation, we
only construct the mapping matrix $A$ for the scattered field,
\begin{equation}
\beta_{p}^{s}=\sum_{m=-p}^{p}A_{p,m}\bar{\alpha}_{m},
\end{equation}
where
\begin{equation}
A_{p,m}=\left(  H_{m-p}({k}\tilde{\rho}_{im})e^{i(m-p)\tilde{\theta}_{im}%
}+\int_{0}^{\infty}H_{m-p}({k}\hat{\tilde{\rho}}_{im})e^{i(m-p)\hat
{\tilde{\theta}}_{im}}\mu(s)ds\right)  . \label{m2l}%
\end{equation}

Notice that the integrand in Eq. (\ref{m2l}) is highly oscillatory for large
$s$ and its numerical computation usually requires special treatment, for
instance, by choosing different integration contours. In the following, we use
the Sommerfeld representations of $H_{n}({k}\rho)e^{in\theta}$ in
Eq.~(\ref{eq:higherordersommer}), but separate it to the propagating and
evanescent terms as in Eq.~(\ref{eq:evanprop}), and reformulate $H_{m-p}({{k}%
}\tilde{\rho}_{im})e^{i(m-p)\tilde{\theta}_{im}}$ using the \textit{plane
wave} representation
\begin{align}
H_{m-p}({k}\tilde{\rho}_{im})e^{i(m-p)\tilde{\theta}_{im}}    =&\frac{i^{m-p}%
}{\pi}\int_{0}^{\pi}e^{i{k}(y\sin{\tau} - x\cos{\tau})}e^{-i(m-p)\theta}%
d\tau\nonumber\\
&  + \frac{(-i)^{m-p}}{i\pi}\int_{0}^{\infty}\frac{e^{-ty}}{\sqrt{t^{2}%
+{k}^{2}}}K(t)dt,
\end{align}
where $(x,y)$ are the Cartesian coordinates of $(\tilde{\rho}_{im}, \tilde{\theta}_{im})$ and
\[
K(t)=e^{i\sqrt{t^{2}+{k}^{2}}x}\left(  \frac{\sqrt{t^{2}+{{k}}^{2}}-t}{{k}%
}\right)  ^{m-p}+e^{-i\sqrt{t^{2}+{k}^{2}}x}\left(  \frac{-\sqrt{t^{2}+{k}%
^{2}}-t}{{k}}\right)  ^{m-p}%
\]
for $y>0$. This plane wave representation was also used to diagonalize the M2L
translation operator in the new version of the low frequency FMM for the
free-space Green's function in Ref. \cite{greengard1998accelerating}. We skip
the similar formula for $H_{m-p}({k}\hat{\tilde{\rho}}%
_{im})e^{i(m-p)\hat{\tilde{\theta}}_{im}}$, and present the translation matrix
explicitly as
\begin{align}
A_{p,m}=  &  \frac{i^{m-p}}{\pi}\int\limits_{0}^{\pi}e^{i{k}(y\sin{\tau} - x\cos{\tau})}e^{-i(m-p)\theta}\left(  \frac{{k}\sin(\tau)-\alpha}{{k}\sin
(\tau)+\alpha}\right)  d\tau\nonumber\\
&  +\frac{(-i)^{m-p}}{i\pi}\int\limits_{0}^{\infty}\frac{e^{-ty}}{\sqrt
{t^{2}+{k}^{2}}}\left(  e^{i\sqrt{t^{2}+{k}^{2}}x}\left(\frac{\sqrt{t^{2}+{k}^{2}%
}-t}{{k}}\right)^{m-p}\right. \nonumber\\
&  \left.  +e^{-i\sqrt{t^{2}+{k}^{2}}x}\left(\frac{-\sqrt{t^{2}+{k}^{2}}-t}{{k}%
}\right)^{m-p}\right)  \left(  \frac{t+i\alpha}{t-i\alpha}\right)  dt
\label{eq:m2loperator}%
\end{align}
after integrating the $s$ variable analytically. In the numerical
evaluation, as the integral of the propagating term is over a finite interval,
high order Gauss quadrature can be applied and for the evanescent term, the
generalized Laguerre quadrature with weight function $t^{n}e^{-t}$ is used.

The translation matrix $A$ in Eq.~(\ref{eq:m2loperator}) has several special
features. Unlike in the classical FMM algorithms, it depends on the $x$ and
$y$ and is therefore spatially variant. However, for the fixed $x$ and $y$, it
is not a two variable function of $m$ and $p$ and only depends on $m-p$. In
the numerical implementation, the matrix can be either computed on-the-fly
using high order Gauss and Laguerre quadratures, or precomputed and stored. We
can estimate the required storage in the algorithm as follows: the translation
operators $A_{p,m}$ are needed for all levels of the tree. For a fixed box,
translation matrix consist of $4p$ complex values (as it is only a function of
$m-p$) and there are a total of no more than $7\cdot7 = 49$ surrounding boxes
representing the well-separated ``receiving" boxes. For the two-layered media
case, the matrix also depends on the $y$-coordinate of the center of the box
as the translation operator takes different values as their images change.
Thus, we can conclude that at tree level $l$, there are $2^{l}$ different
values of $y$-coordinates, and for each $y$-coordinate $49$ possible
well-separated boxes that requires $4p$ complex values. Therefore, the total
required storage for a system with $L$-levels is approximately $(2^{L+1}
\cdot49 \cdot4p )\cdot16 $ bytes, which is very small compared with the
required storage for different expansions.

\subsection{Accelerated evaluation of local direct interactions}

We consider the submatrix representing the \textquotedblleft local direct"
interactions in this section. For a source box with $N_{b}$ particles located
at $\{\mathbf{x}_{j}\}_{j=1}^{N_{b}}$, its domain Green's function
contribution to a target point $\mathbf{x}$ in a neighboring box is defined as
(using the complex image representation)
\begin{align}
u^{d}(\mathbf{x})  &  =\frac{i}{4}\sum_{j=1}^{N_{b}}q_{j}H_{0}({k}%
|\mathbf{x}-\mathbf{x}_{j}|)\nonumber\\
&  +\frac{i}{4}\sum_{j=1}^{N_{b}}q_{j}\left(  H_{0}^{(1)}({k}|\mathbf{x}%
-\mathbf{x}_{j}^{im}|)+\int_{0}^{\infty}H_{0}^{(1)}({k}|\mathbf{x}%
-(\mathbf{x}_{j}^{im}-s\hat{y})|)\mu(s)ds\right). \label{eq:direct}%
\end{align}
This formula shows the entries in one row of the submatrix. Further
compression of the submatrix is usually impossible as it is not low-rank.
However, it is still possible to take advantage of the compressed scattered
field representations of the domain Green's function given in
Eq.~(\ref{eq:greenimage}) (complex image representation) or
Eq.~(\ref{eq:compressed}) (Sommerfeld integral representation), so the entries
in the submatrices can be evaluated more efficiently.

\begin{figure}[t]
\centering  \includegraphics[width=3.0in]{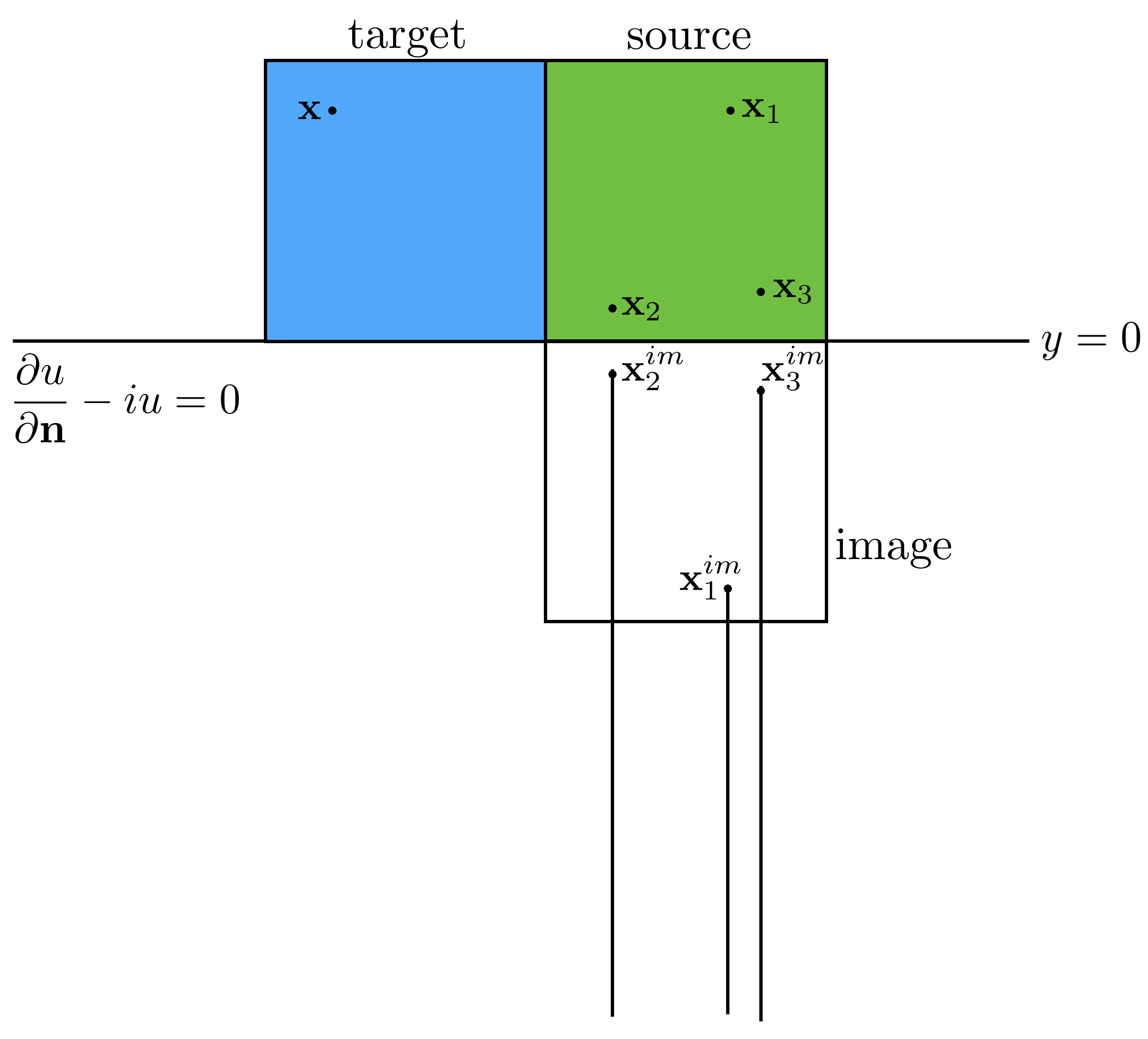} \caption{Images
are separated to near- and far-field by choosing appropriate $C$.}
\label{fig:neighbor}
\end{figure}In Fig.~\ref{fig:neighbor}, we show a source box sitting next to a
target box. In this figure, notice that most of the line-images are
well-separated from the target box. We can therefore choose an appropriate
constant $C$ and cut the line-images into two parts: those that are
well-separated from the target box and those that are not. The evaluation of
Eq.~(\ref{eq:direct}) can be divided as $u^{d}(\mathbf{x})=I+II$, where
\begin{align}
I= &  \frac{i}{4}\sum_{j=1}^{N_{b}}q_{j}H_{0}({k}|\mathbf{x}-\mathbf{x}%
_{j}|)\\
&  +q_{j}\left(  H_{0}^{(1)}({k}|\mathbf{x}-\mathbf{x}_{j}^{im}|)+\int_{0}%
^{C}H_{0}^{(1)}({k}|\mathbf{x}-(\mathbf{x}_{j}^{im}-s\hat{y})|)\mu
(s)ds\right)  ,\nonumber\\
II =&\frac{i}{4}\sum_{j=1}^{N_{b}}q_{j}\int_{C}^{\infty}H_{0}^{(1)}%
({k}|\mathbf{x}-(\mathbf{x}_{j}^{im}-s\hat{y})|)\mu(s)ds. \label{eq:part12}%
\end{align}

The first summation $I$ is computed directly using high order quadrature for
the finite size integral. For the second summation $II$, because
$\mathbf{x}_{j}^{im}-s\hat{y}$ is far away from the target point, the
computation can be accelerated using the available source box multipole
expansion as follows.
\begin{align}
II  &  =\frac{i}{4}\sum_{j=1}^{N_{b}}q_{j}\int_{C}^{\infty}H_{0}^{(1)}%
({k}|\mathbf{x}-(\mathbf{x}_{j}^{im}-s\hat{y})|)\mu(s)ds\nonumber\\
&  =\frac{i}{4}\int_{C}^{\infty}\sum_{m=-P}^{P}\bar{\alpha}_{m}H_{m}%
({k}|\mathbf{x}-(\mathbf{x}_{c}^{im}-s\hat{y})|)e^{im\hat{\theta}_{im}}%
\mu(s)ds\nonumber\\
&  =\frac{i}{4}\int_{C}^{\infty}\sum_{m=-P}^{P}\bar{\alpha}_{m}\sum
_{n=-\infty}^{\infty}H_{m-n}({k}\tilde{\rho}_{im})e^{i(m-n)\hat{\tilde{\theta
}}_{im}}J_{n}({k}|\mathbf{x}-\mathbf{x}_{c}^{l}|)e^{in\theta}\mu
(s)ds\nonumber\\
&  =\frac{i}{4}\sum_{n=-\infty}^{\infty}\left(  \sum_{m=-P}^{P}\bar{\alpha
}_{m}\int_{C}^{\infty}H_{m-n}({k}\tilde{\rho}_{im})e^{i(m-n)\hat{\tilde
{\theta}}_{im}}\mu(s)ds\right)  J_{n}({k}|\mathbf{x}-\mathbf{x}_{c}%
^{l}|)e^{in\theta}\nonumber\\
&  \approx\frac{i}{4}\sum_{n=-P}^{P}L_{n}J_{n}({k}|\mathbf{x}-\mathbf{x}%
_{c}^{l}|)e^{in\theta},
\end{align}
where
\[
L_{n}=\sum_{m=-P}^{P}\bar{\alpha}_{m}\int_{C}^{\infty}H_{m-n}({k}\tilde{\rho
}_{im})e^{i(m-n)\hat{\tilde{\theta}}_{im}}\mu(s)ds=\sum_{m=-k}^{k}\bar{\alpha
}_{m}B_{m,k}%
\]
and the translation matrix is given by
\begin{equation}
B_{m,k}=\int_{C}^{\infty}H_{m-n}({k}\tilde{\rho}_{im})e^{i(m-n)\hat
{\tilde{\theta}}_{im}}\mu(s)ds,
\end{equation}
which can be evaluated efficiently using the corresponding Sommerfeld integral representation.

In the hierarchical tree structure, most boxes are well-separated from the
interface $y=0$. This implies that $C=0$ for most direct interactions of the
source and target boxes, and the separation can be simplified as
\begin{align}
I  &  =\frac{i}{4}\sum_{j=1}^{N_{b}}q_{j}H_{0}({k}|\mathbf{x}-\mathbf{x}%
_{j}|),\nonumber\\
II  &  =\frac{i}{4}\sum_{j=1}^{N_{b}}q_{j}\left(  H_{0}^{(1)}({k}%
|\mathbf{x}-\mathbf{x}_{j}^{im}|)+\int_{0}^{\infty}H_{0}^{(1)}({k}%
|\mathbf{x}-(\mathbf{x}_{j}^{im}-s\hat{y})|)\mu(s)ds\right)  ,
\end{align}
where both the point- and line-image contributions belong to $II$. In this
case, the corresponding translation operator becomes the same as in
Eq.~(\ref{eq:m2loperator}). In the numerical simulation, all the translation
matrices can be either precomputed or computed on-the-fly using high order
quadrature for the Sommerfeld integral representation.

The selected compression schemes and translations allow easy adaptation of
existing fast multipole algorithms for computing the domain Green's function
interactions of the 2-D half-space Helmholtz equation with impedance boundary
condition. Our solver is based on the wide-band FMM
\cite{cho2010wideband,cho2012revision}. We present the pseudo-code of our
algorithm in {\bf Algorithm 1}.

\begin{algorithm}
\caption{Heterogeneous 2-D FMM for Two-layered Media with Impedance Boundary Conditions}
\begin{algorithmic}
\STATE{\begin{center}{\bf Step 1: Initialization}\end{center}}
\STATE {\bf Generate} an adaptive hierarchical tree structure and precompute tables.
\STATE{{\bf Comment [} $L$ denotes the maximum refinement level in the adaptive tree
	determined by a prescribed number $s$ representing the maximum allowed number of
	particles in a childless box.
{\bf ]}}
\STATE{\begin{center}{\bf Step 2: Upward Pass}\end{center}}
\FOR{$l = L, \cdots 0$}
\FOR{ all boxex $j$ on level $l$}
\IF {$j$ is a leaf node}
\STATE{{\bf form} the {\it free-space} multipole expansion using Eq.~(\ref{multipole_free}).}
\ELSE
	\STATE{{\bf form} the {\it free-space} multipole expansion by merging children's expansions using the {\it free-space} ``multipole-to-multipole" translation operator.}
\ENDIF
\ENDFOR
\ENDFOR
\STATE{{\bf Cost [} All operations in this step are the same as those in the {\it free-space} FMM. {\bf ]} }
\STATE{\begin{center}{\bf Step 3: Downward Pass}\end{center}}
\FOR{$l=1, \cdots, L$}
	\FOR{all boxes $j$ on level $l$}
		\STATE {{\bf shift} the local expansion of $j$'s parent to $j$ itself using the {\it free-space} ``local-to-local" translation operator.}
	\STATE{{\bf collect} interaction list contribution using the precomputed table and the ``multipole-to-local" translation operator in Eq.~(\ref{m2l}).}
	\ENDFOR	
\ENDFOR
\STATE{\bf Cost [} Using the precomputed table, the cost is expected to be the same as in the {\it free-space} FMM. Overhead
operations are required when tables are computed on-the-fly. {\bf ]}
\STATE{\begin{center}{\bf Step 4: Evaluate Local Expansions}\end{center}}
\FOR{each leaf node (childless box)}
\STATE{{\bf collect} part II in Eq.~(\ref{eq:part12}) from neighboring (including self) boxes.}
\STATE{{\bf evaluate} the local expansion at each particle location.}
\ENDFOR
\STATE{{\bf Comment [} At this point, for each target point, its far field contribution (including those
from well-separated images) has been computed.  {\bf ]} }
\noindent
\STATE{{\bf Cost [} Compared with the {\it free-space} FMM, additional translations are required to translate
the multipole expansions of images to local expansions. The heterogeneous translation operators can
be computed on-the-fly or precomputed. The amount of work is constant for each leaf node.
{\bf ]}}
\STATE{\begin{center}{\bf Step 5: Local Direct Interactions}\end{center}}
\FOR{$i=1, \cdots, N$}
\STATE {{\bf compute} Part I in Eq.~(\ref{eq:part12}) of target point $i$ with original and image  sources in the neighboring boxes.}
\ENDFOR
\STATE{{\bf Cost [} When the computational domain is well-separated from the boundary $y=0$, this step only
involves the evaluation of the {\it free-space} kernel and the cost is the same as the {\it free-space}
FMM. When the computational domain is close to the boundary $y=0$, a constant number of additional
operations are required for each $i$ in a very small subset of the particles to evaluate the
near-field point- and line-image contributions from Part I in Eq.~(\ref{eq:part12}).
{\bf ]}}
\end{algorithmic}
\end{algorithm}

In our current implementation, all the tables are precomputed using
Mathematica requesting more than 20 digits accuracy. Compared with the
original free-space FMM algorithm, the domain Green's function FMM only
requires a small portion of additional cost, as demonstrated in the next section.

\section{Numerical Results}

\label{sec:numerical_results}

\begin{figure}[t]
\centering
\includegraphics[width=4.2in]{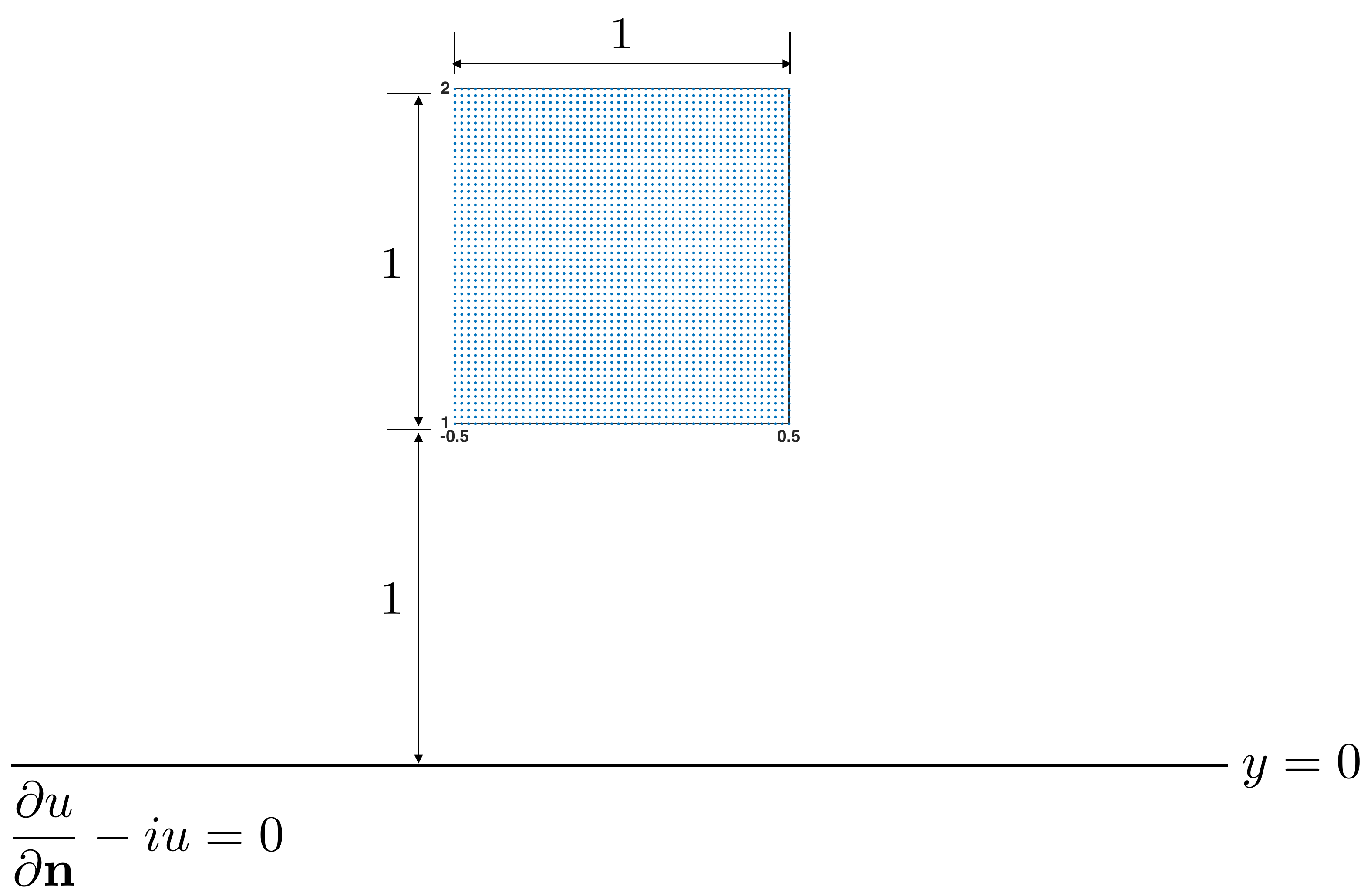}\caption{Uniform distribution in a unit
square on top of half-space.}%
\label{example1}%
\end{figure}

We present some preliminary numerical results in this section to demonstrate
the performance of the new heterogeneous FMM algorithm for the two-layered
media with the interface placed at $y = 0$, and set $\alpha= 1$ in the
impedance boundary condition. We assume the source and target points are the
same set of $N$ particles located in a unit box centered at $(0, 1.5)$ as
shown in Fig. \ref{example1}. The numerical simulations are performed on a
desktop with 3.7 GHz Xeon E5 processor and 32GB RAM using the gcc compiler
version 4.9.3. All the required translation tables are precomputed using
Mathematica. 

\begin{figure}[ptb]
\caption{CPU time (seconds) for different $N$ using $p = 39$ and ${k} = 0.1$}
\centering  \includegraphics[width=4.5in]{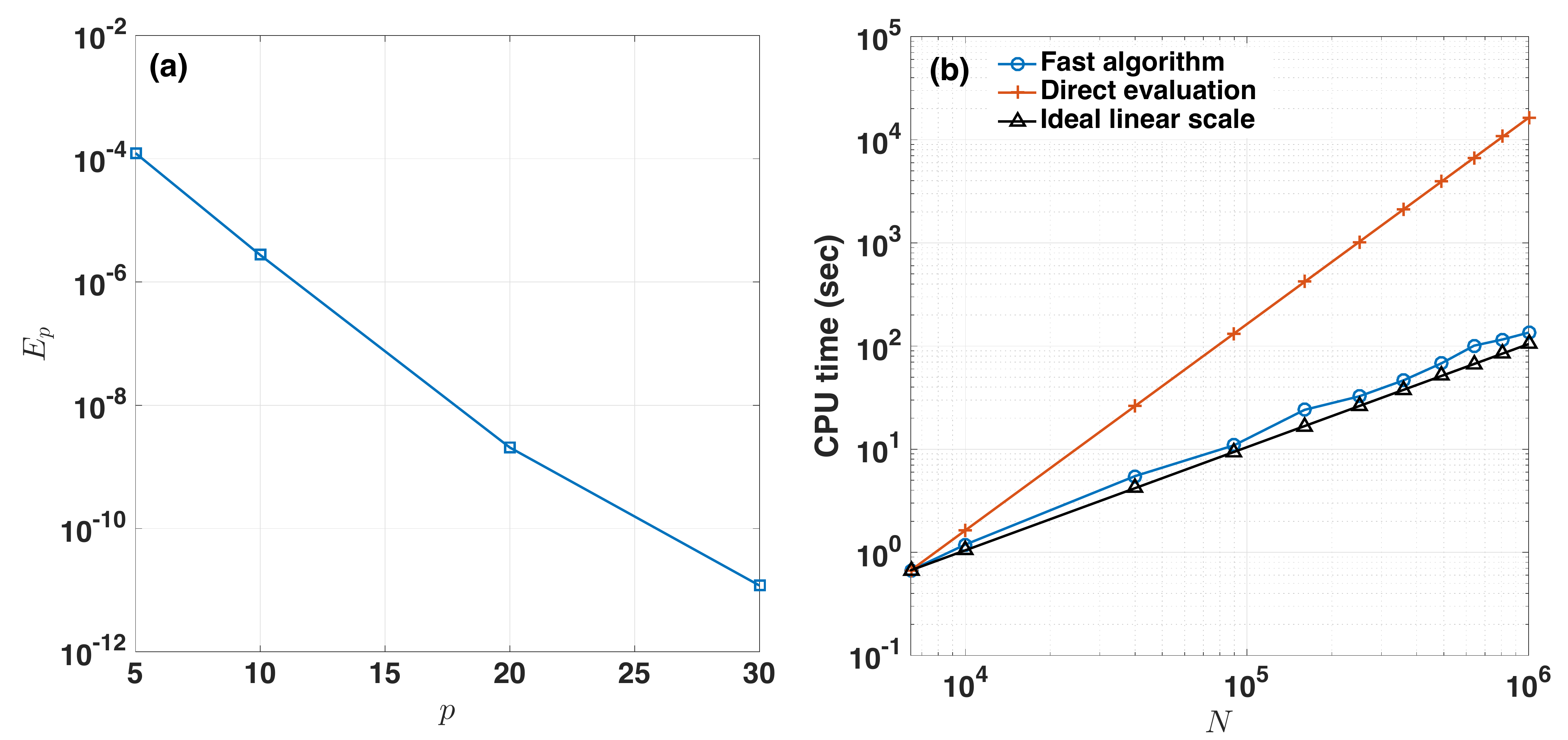} \caption{(a)
Convergence and (b) linear CPU time scaling for the impedance half-space
problem with ${k} = 0.1$.}%
\label{CPU_Error_plot}%
\end{figure}

As the analytical solution is not available for this problem, we first check
the algorithm accuracy by studying how the errors change as a function of the
number of expansion terms $p$. We consider the example with $N=100\times100$
particles uniformly distributed in the box. A reference solution is computed
using $p=39$ (which should provide results with approximately $12$-digit
accuracy) and by setting $L=3$ in the hierarchical tree structure. The new
heterogeneous FMM algorithm took about $1.19$ seconds to derive the reference
solutions. The accuracy results are presented in Table \ref{accuracy_table}
for $p=5,10,20,30$ and wave numbers ${k}=0.1$ and ${k}=1$, respectively, where
the error $E_{p}$ using $p$ terms in the expansion is defined as
\begin{equation}
E_{p}=\left(  \frac{\sum_{j=1}^{M}|u_{39}(\mathbf{x}_{j})-u_{p}(\mathbf{x}%
_{j})|^{2}}{\sum_{j=1}^{M}|u_{39}(\mathbf{x}_{j})|^{2}}\right)  ^{\frac{1}{2}%
},\text{ \ \ }M=10,000.\label{accuracy}%
\end{equation}
In Fig. \ref{CPU_Error_plot}(a), we plot how the error $E_{p}$ decays as a
function of $p$ for ${k}=0.1$. We see that the error dependency on the number
of terms $p$ to compress the domain Green's function is similar to that in
existing FMM analysis for the free-space kernels. \begin{table}[t]
\caption{Accuracy results with different expansion terms for ${k}=0.1$ and
${k}=1$. Reference solution is computed with $p=39$.}%
\label{accuracy_table}
\centering
\begin{tabular}
[c]{|c|c|c|}\hline
$p$ & Error for ${k} = 0.1$ & Error for ${k} = 1$\\\hline
$E_{5}$ & $1.23\times10^{-4}$ & $1.43 \times10^{-4}$\\\hline
$E_{10}$ & $2.73\times10^{-6}$ & $3.81 \times10^{-6}$\\\hline
$E_{20}$ & $2.06\times10^{-9}$ & $2.85 \times10^{-9}$\\\hline
$E_{30}$ & $1.19\times10^{-11}$ & $1.65 \times10^{-11}$\\\hline
\end{tabular}
\end{table}

We demonstrate the algorithm efficiency by presenting the CPU times in Table
\ref{timing_table} for different numbers of source/target points $N$ from
$100$ to $1,000,000$ for ${k} = 0.1$. A $\log$-$\log$ plot of the CPU time is
also presented in Fig. \ref{CPU_Error_plot}(b), which clearly shows the linear
scaling of the new heterogeneous FMM algorithm. For comparisons, estimated
results of direct computations (using CPU times for $N=100$ and $6,400$) and
the ideal linear scaling curve are also presented. Similar experiments are
performed for ${k} = 1$ and results are almost identical to that when ${k} =
0.1$ and are therefore omitted in this paper. \begin{table}[t]
\caption{CPU time (seconds) for different $N$ using $p = 39$ and ${k} = 0.1$}%
\label{timing_table}
\centering
\par%
\begin{tabular}
[c]{|c|c|c|c|c|c|c|c|c|c|c|}\hline
$N$ & 100 & 6400 & 10000 & 90000 & 360000 & 640000 & 810000 & 1000000
\\\hline
CPU time & 0.01 & 0.67 & 1.19 & 10.92 & 46.58 & 100.85 & 116.03 & 135.05
\\\hline
\end{tabular}
\end{table}

\section{Heterogeneous FMM for Three Layers}

\label{sec:threelayer_fmm}

In this section, we study a three-layered media setting where the source and
target points are all located in the top layer and their interactions are
described by Eq.~(\ref{eq:multitop}), and present the framework of a
hierarchical algorithm based on the Sommerfeld integral representation. This
framework allows direct generalization to multiple layers with analysis and
algorithm design for more general settings. All the technical details of the
compressions, translation operators, and bookkeeping strategies for different
numbers of layers, boundary conditions, and locations of the source and target
points will be presented in a subsequent paper. \newline

{\noindent\textbf{Compressed Representation:}} We start from the Sommerfeld
integral representation of the domain Green's function for sources and target
points located in the top layer described in Eq.~(\ref{eq:multitop}). Similar
to Eqs.~(\ref{eq:compressed}) and (\ref{eq:uncompressed}) in
Sec.~\ref{sec:lowrank} for the two-layered case, we first derive the
compressed ``multipole" expansion in the form
\begin{align}
u_{1}^{s}(\mathbf{x})\approx &  \int_{-\infty}^{\infty}\underbrace
{\frac{e^{-\sqrt{\lambda^{2}-{k^{2}_{1}}}(y+y_{c})}}{\sqrt{\lambda^{2}-{k_{1}%
}^{2}}}e^{i\lambda(x-x_{c})}}_{\mbox{free-space info}}\underbrace{\left(
\frac{1}{4\pi}\sum_{p=-P}^{P}\bar{\alpha}_{p}(-i)^{p}\left(  \frac
{\lambda-\sqrt{\lambda^{2}-{k^{2}_{1}}}}{{k_{1}}}\right)  ^{p}\right)
}_{\mbox{compressed}}\label{eq:multicompressed}\\
&  \underbrace{ \sigma_{1}(\lambda) }_{\mbox{image info}}d\lambda,\nonumber
\end{align}
where $(x_{c},y_{c})$ represents the center of the source box, and the
coefficients $\bar{\alpha}_{p}$ are either computed directly using the
free-space ``particle-to-multipole" translation described in
Eq.~(\ref{multi_coeff}) for a childless box, or computed using the free-space
``multipole-to-multipole" translation operator for a parent box.

\medskip
{\noindent\textbf{Remark 4:}} We expect that the error analysis for the three
layers, in fact for multiple layers for this matter, should follow in a manner
similar to the case of two layers because the point and line images in the domain Green's
function used to compute the far field will stay far away from the far field
box as illustrated in Fig. 1. We will present the construction of the images for multiple
layers and the mathematical error analysis of compression and translation operators
in the follow-up paper.

Moreover, it should be noted that the error analysis can also be performed without
using the image representation, by observing that the decay rate of the series
expansion is determined by the ratio of $||(x,y)-(x_{c},-y_{c})||$ and the
size of the source box. As the source and target boxes are well-separated and
the distance $||(x,y)-(x_{c},-y_{c})||>||(x,y)-(x_{c},y_{c})||$, we conclude
that the number of terms required in the expansion is no more than that in the free-space.

\medskip
{\noindent\textbf{Multipole-to-local Translation:}} Similar to the two-layered
media, we use the Hankel function based free-space ``local" expansions to
collect the far-field compressed representation. There are several ways to
derive the transformation matrix from the multipole coefficients $\bar{\alpha
}_{p}$ in Eq.~(\ref{eq:multicompressed}) to the local coefficients. For
instance, one can use the identity (see \cite{abramowitz1966handbook})
\[
e^{i k r \cos(\theta)} = \sum_{m=-\infty}^{\infty} i^{m} e^{i m \theta}
J_{m}(k r)
\]
to derive the local expansion coefficients after shifting the multipole
expansion in Eq.~(\ref{eq:multicompressed}) to the new center of the target
box. This approach was also used in the new version of FMM algorithm for the
Helmholtz equation \cite{greengard1998accelerating}. A much simpler approach
is to compare with the two-layered media formula and modify
Eq.~(\ref{eq:m2loperator}) to directly derive the translation matrix as
\begin{align}
A_{p,m}=  &  \frac{i^{m-p}}{\pi}\int\limits_{0}^{\pi}e^{i{k_{1}}(y\sin{\tau} - x\cos{\tau})}e^{-i(m-p)\theta}\left(  \sigma_{1}\left( -k_{1} \cos
(\tau)\right)  \right)  d\tau\nonumber\\
&  +\frac{(-i)^{m-p}}{i\pi}\int\limits_{0}^{\infty}\frac{e^{-ty}}{\sqrt
{t^{2}+{k^{2}_{1}}}}\left(  e^{i\sqrt{t^{2}+{k^{2}_{1}}}x}\left(\frac{\sqrt
{t^{2}+{k^{2}_{1}}}-t}{{k_{1}}}\right)^{m-p}\right. \nonumber\\
&  ~\left.  +e^{-i\sqrt{t^{2}+{k^{2}_{1}}}x}\left(\frac{-\sqrt{t^{2}+{k^{2}_{1}}%
}-t}{{k_{1}}}\right)^{m-p}\right)  \left(  \sigma_{1}\left(  \sqrt{t^{2}+k^{2}_{1}%
}\right)  \right)  dt.
\end{align}
The translation matrix can be either precomputed, or computed on-the-fly using
high order Gauss and Laguerre quadratures. Once the local expansion is
available, the ``local-to-local" translations are the same as the free-space ones.

\medskip

{\noindent\textbf{Accelerating the evaluation of local direct interactions.}}
Similar to the two-layered case, the local direct interactions can be
separated into the free-space interactions and scattered field due to the
images. However, finding the corresponding explicit image representation from
the Sommerfeld integral usually requires a numerical integral transformation
and how to reduce its computational complexity is still a research topic. We
therefore use the following strategies to accelerate the evaluation of local
direct interactions: (a) when the source box is located at least one box size
away from the interface $y=0$, the compressed multipole expansion in
Eq.~(\ref{eq:multicompressed}) for the scattered field is then valid at the
target box, we therefore perform a \textquotedblleft multipole-to-local"
translation using a precomputed table, which will add the scattered field
contribution to the local expansion of the box; (b) otherwise, we compute the
direct interactions using numerical integration schemes for the domain Green's
function.\medskip

{\noindent\textbf{Remark 5:}} Mathematically, the cut-off idea is the same as
separating the term $\sigma_{1}$ in Eq.~(\ref{eq:multitop}) as the sum of two
terms, so that one term contains a factor $e^{-\sqrt{\lambda^{2}-k_{1}^{2}}C}$
which makes the truncated \textquotedblleft multipole expansion" in
Eq.~(\ref{eq:multicompressed}) valid for a prescribed accuracy requirement,
while the other term when transformed into the physical domain, only involved
a finite integral where the \textquotedblleft images" are located in a small
region in the neighborhood of the target box. The numerical separation of
$\sigma_{1}$ may be more difficult than computing the Sommerfeld
representation directly.\medskip

{\noindent\textbf{Algorithm Summary.}} The hierarchical algorithm for the
special three-layered media setting can be summarized as follows. It uses the
same tree structure, free-space compressed multipole coefficient,
multipole-to-multipole, and local-to-local translation operators as the
two-layered algorithm. Also, when the source box is well separated from the
interface $y=0$, its local direct interactions can be accelerated similarly as
in the two-layered media case. There are only two necessary revisions of the
two-layered media algorithm: (a) a different set of tables for the
\textquotedblleft multipole-to-local" translation operators are used; and (b)
when the \textquotedblleft images" in the scattered field are no longer
well-separated from the target boxes when evaluate the local direct
interactions, instead of using a cut-off distance, the domain Green's function
is evaluated directly using some numerical integration schemes. As the
algorithm structure is almost identical to the two-layered case, we neglect
the details in this paper.\medskip

{\noindent\textbf{Remark 6:}} In the numerical experiment, when the
configuration of the source and target points in the top layer is the same as
that in the upper layer of the half-space problem in Fig.~\ref{example1}, the
three-layered media algorithm only differs from two-layered one in the
precomputed tables, and both the accuracy and efficiency are therefore the
same and we neglect the details.

\section{Conclusion and Future Work}

\label{sec:conclusion}

We present a heterogeneous FMM for the efficient calculation of
the discretized integral operator for the Helmholtz equation in two and three
layers media with details on the case of two-layered media with impedance boundary
conditions. The two-layered media setting
allows the use of both the \textit{complex line image} and \textit{Sommerfeld
integral} representations to compress the domain Green's function and to
derive the translation operators analytically. Instead of compressing the
interaction matrix directly, the \textit{complex line image} representation
intuitively reveals how a ``transformed" matrix can be compressed through a
procedure that only involves the free-space Green's function, and provides
rigorous error bounds by using existing free-space FMM results.

Unlike the fast direct solvers, the compression is performed analytically on a
\textquotedblleft transformed" matrix which allows the easy adaptation of
existing free-space FMM packages. Also disimilar to the classical FMM, the
\textquotedblleft multipole-to-local" translation operators are spatially
variant, thus the translation operators in the FMM becomes heterogeneous.
Numerical experiments show that the new hierarchical algorithm provides
significant improvement over existing hybrid methods \cite{o2014efficient} in
two-layered media settings. By connecting the complex line image and
Sommerfeld representations in the two-layered media algorithm, we demonstrate
how the compressions and translations can be performed directly on the
Sommerfeld representation for a particular setting in the three-layered media
for the case when all the source and targets points are located in the top layer.

This paper focuses on the intuitions through the two-layered and three-layered
media setting. In a subsequent paper, we will present the H-FMM for multiple
layered media including constructions and error analysis of  compressions,
translations as well as bookkeeping strategies, for both the scalar Helmholtz
equations in acoustic studies and the multi-layered media dyadic Green's
function for the Maxwell's equations \cite{cho2016efficient}.

Finally, it is interesting to compare the analysis based compressions with those using purely
numerical linear algebra techniques as in the fast direct solvers, to
understand how the efficiencies of both compressions can be further improved.
Research along these directions will be explored in the future.


\bibliographystyle{siamplain}
\bibliography{HFMM}

\end{document}

%% file: ex_shared.tex
\usepackage{lipsum}
\usepackage{amsfonts}
\usepackage{graphicx}
\usepackage{epstopdf}
\usepackage{algorithmic}
\ifpdf
  \DeclareGraphicsExtensions{.eps,.pdf,.png,.jpg}
\else
  \DeclareGraphicsExtensions{.eps}
\fi

\numberwithin{theorem}{section}

\newcommand{\TheTitle}{A Heterogeneous FMM for 2-D Layered Media Helmholtz Equation I: Two \& Three Layers Cases}
\newcommand{\TheAuthors}{M.H. Cho, J. Huang, D. Chen, and W. Cai}

\headers{A Heterogeneous FMM for 2-D Layered Media}{\TheAuthors}

\title{{\TheTitle}\thanks{Submitted to the editors DATE.
\funding{This work was supported by a grant from the Simons Foundation (\#404499, Min
Hyung Cho). J. Huang was supported by NSF (Grant No. DMS-1217080), and W. Cai
is supported by Army Research Office (Grant No. W911NF-14-1-0297) and NSF
(Grant No. DMS-1619713).}}}

\author{
  Min Hyung Cho\thanks{Department of Mathematical Sciences, University of Massachusetts Lowell, Lowell, MA 01854,
    (\email{minhyung\_cho@uml.edu})}
  \and
  Jingfang Huang \thanks{Department of Mathematics, University of North
    Carolina at Chapel Hill, Chapel Hill, NC 27599-3250. (\email{huang@email.unc.edu}, \email{dangxing@live.unc.edu})}
    \and
   Dangxing Chen \footnotemark[3]
  \and
  Wei Cai\thanks{Corresponding author, Beijing Computational Science Research Center, Beijing, China, (\email{wcai@csrc.ac.cn}); Department of Mathematics and Statistics, University of North Carolina at Charlotte, Charlotte, NC 28223 (\email{wcai@uncc.edu}).}
}

\usepackage{amsopn}
